\documentclass[12pt, reqno]{amsart}
\usepackage{amsfonts}
\usepackage{amssymb}
\usepackage{amsmath}

\newtheorem{theorem}{Theorem}
\theoremstyle{plain}

\newtheorem{remark}{Remark}

\begin{document}
\title[] {Lebesque-Radon-Nikodym theorem with respect to fermionic $p$-adic invariant measure on $\Bbb Z_p$}
\author{Taekyun Kim}
\address{Taekyun Kim. Division of General Education-Mathematics \\
Kwangwoon University, Seoul 139-701, Republic of Korea  \\}
\email{tkkim@kw.ac.kr}
 \maketitle

{\footnotesize {\bf Abstract} \hspace{1mm} {In this paper we derive
the analogue of the Lebesque-Radon-Nikodym theorem with respect to
fermionic $p$-adic invariant measure on $\Bbb Z_p$.}



\section{Introduction}
Let $p$ be a fixed odd prime number. Throughout this paper, the
symbol $\Bbb Z$, $\Bbb Z_p$, $\Bbb Q_p$, and $\Bbb C_p$ denote the
ring of rational integers, the ring of $p$-adic integers, the field
of $p$-adic rational numbers, and the completion of algebraic
closure of $\Bbb Q_p$, respectively. Let $\nu_p $ be the normalized
exponential valuation of $\Bbb C_p$ with $|p|_p=
p^{-\nu_{p}(p)}=p^{-1}$. In this paper, we assume that  $ q \in \Bbb
C_p$ with $|q-1| < p^{- \frac{1}{p-1}}$ as an indeterminate. We use
the notations
$$ [x]_q = [x:q]=\frac{1- q^x }{1-q},~~ [x]_{-q} =
\frac{1-(-q)^x }{1+q},$$  for all $x \in \Bbb Z_p$. For any positive
integer $N$, let $a+ p^N \Bbb Z_p = \{ x\in X | \,  x \equiv a
\pmod{p^N}\}$, where $a\in \Bbb Z$ lies in $0\leq a <  p^N$. The
$p$-adic $q$-invariant distribution $\mu_q$ on $\Bbb Z_p$ is defined
by
$$\mu_q( a+ p^N \Bbb Z_p) = \frac{q^a}{[p^N ]_q}, \quad \text{(see~~[4-8])}.$$

We say that $f$ is a uniformly differentiable function at a point $a
\in \Bbb Z_p$ and denote this property by $f \in UD(\Bbb Z_p)$, if
the difference quotients $F_f (x,y)= \frac{f(x)-f(y)}{x-y}$ have a
limit $l=f'(a)$ as $(x,y) \rightarrow (a,a)$. For $f \in UD(\Bbb
Z_p)$, the $p$-adic invariant $q$-integral is defined as
\begin{eqnarray*} I_{q} (f)=\int_{\Bbb
Z_p} f(x) d \mu_{q} (x)= \lim_{N \rightarrow \infty}
\frac{1}{[p^N]_{q}} \sum_{x=0}^{p^N -1} f(x) q^x.
\end{eqnarray*}
The fermionic $p$-adic $q$-measures on $\Bbb Z_p$ are defined as
$$\mu_{-q}( a+d p^N \Bbb Z_p) = \frac{(-q)^a}{[d p^N ]_{-q}},$$
and the fermionic $p$-adic invariant $q$-integral on $\Bbb Z_p$ is
defined as
\begin{eqnarray*} I_{-q} (f)=\int_{\Bbb
Z_p} f(x) d \mu_{-q} (x)= \lim_{N \rightarrow \infty}
\frac{1}{[p^N]_{-q}} \sum_{x=0}^{p^N -1} f(x) (-q)^x,
\end{eqnarray*}
for $f \in UD(\Bbb Z_p)$. The fermionic $p$-adic invariant integral
on $\Bbb Z_p$ is defined as
$$I_{-1} (f)=\lim_{q \rightarrow 1} \int_{\Bbb Z_p} f(x) d \mu_{-q} (x)=\int_{\Bbb Z_p} f(x) d \mu_{-1} (x), \text{ (see [6])}. $$
 Let $C(\Bbb Z_p, \Bbb C_p)$ be the space of continuous function on
$\Bbb Z_p$ with values in $\Bbb C_p$, provided with norm
$||f||_\infty=\underset{x \in \Bbb Z_p}{\sup} |f(x)|$. The
difference quotient $\bigtriangleup_1 f$ of $f$ is the function of
two variables given by
$$\bigtriangleup_1 f(m, x) = \frac{f(x+m)-f(x)}{m}, $$
for all $x \in \Bbb Z_p$ and $m \in \Bbb Z_p$ with $m \neq 0$. A
function $f:\Bbb Z_p \longrightarrow \Bbb C_p$ is said to be a
Lipschitz function if there exists a constant $M>0$  such that for
all $m \in \Bbb Z_p \setminus \{0 \}$ and $x \in \Bbb Z_p$,
$$|\bigtriangleup_1 f(m, x)| \le M.$$
Here $M$ is called the Lipschitz constant of $f$. The $\Bbb
C_p$-linear space consisting of all Lipschitz function is denoted by
Lip$(\Bbb Z_p, \Bbb C_p)$. This space is a Banach space with respect
to the norm $||f||_1 = ||f||_\infty \vee ||\bigtriangleup_1
f||_\infty$, (see [8]). Various proofs of the Radon-Nikodym theorem
can be found in many books on measure theory, analysis, or
probability theory. Usually they use the Hahn decomposition theorem
for signed measures, the Riesz representation theorem for
functionals on Hilbert spaces, or martingale theory, ( see [9, 10,
11, 12]). Recently, several authors have studied related to
Radon-Nikodym theorems ( see [1, 2, 3]).
 In previous paper [7], author have studied the analogue of the
Lebesque-Radon-Nikodym theorem with respect to $p$-adic
$q$-invariant distribution on $\Bbb Z_p$. The purpose of this paper
is to derive  the Lebesque-Radon-Nikodym's type theorem with respect
to fermionic $p$-adic invariant measures on $\Bbb Z_p$.

\section{Analogue of the
Lebesque-Radon-Nikodym theorem with respect to fermionic $p$-adic
invariant integral on $\Bbb Z_p$}

Let $f \in UD(\Bbb Z_p)$. For any positive integers $a$ and $n$ with
$a<p^n$, define
\begin{eqnarray}
\mu_{f, -1}(a+p^n \Bbb Z_p )=\int_{a+p^n \Bbb Z_p} f(x)d \mu_{-1}
(x),
\end{eqnarray}
where the integral is the fermionic $p$-adic invariant integral on
$\Bbb Z_p$.

It is easy to see that
\begin{eqnarray}
\mu_{f, -1}(a+p^n \Bbb Z_p )&=&  \lim_{m \rightarrow \infty} \sum_{x=0}^{p^m -1} f(a+p^n x)(-1)^{a+p^n x}  \\
&=&  \lim_{m \rightarrow \infty}(-1)^a \sum_{x=0}^{p^{m-n} -1}
f(a+p^n x)(-1)^{x}.\notag
\end{eqnarray}
Thus it follows that
\begin{eqnarray}\mu_{f, -1}(a+p^n \Bbb Z_p )= (-1)^a \int_{\Bbb Z_p} f(a+p^n x) d\mu_{-1}(x).\end{eqnarray}

By (3), we have that for $f, \, g \in UD(\Bbb Z_p, \Bbb C_p)$,
\begin{eqnarray}\mu_{\alpha f + \beta g , \, -1}=\alpha \mu_{f, \,
-1}+\beta \mu_{g, \, -1}, \end{eqnarray} where $\alpha, \, \beta$
are constants. We also see that
\begin{eqnarray} |\mu_{ f, \, -1}(a+p^n \Bbb Z_p)| \le \|f\|_{\infty}. \end{eqnarray}

Here we recall the definition of the strongly fermionic $p$-adic
invariant measure on $\Bbb Z_p$. If $\mu_{-1}$ is satisfied the
following equation
$$| \mu_{-1} (a+p^n  \Bbb Z_p)- \mu_{-1} (a+p^{n+1} \Bbb Z_p) | \le \delta_n ,$$
where $\delta_n \rightarrow 0$ as $n \rightarrow \infty$ and
$\delta_n$ is independent of $a$, then $\mu_{-1}$ is called the
weakly fermionic invariant measure on $\Bbb Z_p$. If $\delta_n$ is
replaced by $cp^{-n}$ with a constant $c$, then $\mu_{-1}$ is called
strongly fermionic $p$-adic invariant measure on $\Bbb Z_p$.

Let $P(x) \in \Bbb C_p [[x]]$ be an arbitrary polynomial. Now we
claim that $\mu_{P, -1}$ is a strongly fermionic $p$-adic invariant
measure on $\Bbb Z_p$. Without a loss of generality, it is enough to
prove the statement for $P(x)=x^k$.

Let $a$ be an integer with $0 \le a < p^n$. Then we see that
\begin{eqnarray}
\mu_{P, \, -1}(a+p^n \Bbb Z_p ) =  \lim_{m \rightarrow \infty}(-1)^a
\sum_{i=0}^{p^{m-n} -1} (a+i p^n)^k (-1)^i.
\end{eqnarray} It is easy to see that
\begin{eqnarray}(a+i p^n)^k=\sum_{l=0}^{k} a^{k-l} \binom{k}{l}(ip^n)^l = a^k + \binom{k}{1}a^{k-1}p^n i+ \cdots + p^{n^k} i^k.\end{eqnarray}
By (6) and (7), we have
\begin{eqnarray}
& &\mu_{P, \, -1}(a+p^n \Bbb Z_p )\\ & & =  \lim_{m \rightarrow
\infty}(-1)^a \{ a^k + k a^{k-1}p^n \sum_{i=0}^{p^{m-n}-1} i (-1)^i
+ \cdots + p^{n^k}\sum_{i=0}^{p^{m-n}-1}i^k (-1)^i \}.\notag
\end{eqnarray}

From the definition of the Euler numbers $E_n$, we note that
\begin{eqnarray}
E_n =  \int_{\Bbb Z_p} x^n d \mu_{-1} (x) =\lim_{m \rightarrow
\infty}\sum_{i=0}^{p^m-1}i^n (-1)^i,
\end{eqnarray}
(see [6]). By (8) and (9), we obtain that
\begin{eqnarray*}
\mu_{P, -1}(a+p^n \Bbb Z_p ) &\equiv& (-1)^a a^k E_0 \qquad \,
(\text{mod} \, p^n)\\
&\equiv& (-1)^a P(a) E_0 \quad \,(\text{mod} \, p^n).
\end{eqnarray*}

Let $x$ be an arbitrary in $\Bbb Z_p$. Let $x \equiv x_n \,\,
(\text{mod} \,  p^n)$ and $x \equiv x_{n+1} \,\, (\text{mod} \,
p^{n+1})$, where $x_n$ and $x_{n+1}$ are positive integers such that
$0 \le x_n < p^n$ and $0 \le x_{n+1} < p^{n+1}$. Then we have
\begin{eqnarray*}
|\mu_{P, -1}(x+p^n \Bbb Z_p )-\mu_{P, -1}(x+p^{n+1} \Bbb Z_p)| \le C
p^{-n},
\end{eqnarray*}
where $C$ is some constant and $n \gg 0$.

Let
\begin{eqnarray*}
f_{\mu_{P, -1}}(a) =  \lim_{n \rightarrow \infty} \mu_{P, -1}(a+p^n
\Bbb Z_p).
\end{eqnarray*}
Then we have from (8) and (9) that
\begin{eqnarray}
f_{\mu_{P, -1}}(a)
&=& (-1)^a \lim_{n \rightarrow \infty} \{ a^k + ka^{k-1}p^nE_1 + \cdots + p^{n^k} E_k \} \notag \\
&=& (-1)^a a^k = (-1)^a P(a).
\end{eqnarray}
Since $f_{\mu_{P, -1}} (x)$ is continuous, it follows that for all
$x \in \Bbb Z_p$,
\begin{eqnarray}f_{\mu_{P, -1}} (x)=(-1)^x P(x).\end{eqnarray}

For $g \in UD(\Bbb Z_p, \Bbb C_p)$, we see that
\begin{eqnarray}
\int_{\Bbb Z_p} g(x) d \mu_{P, -1} (x) &=&  \lim_{n \rightarrow \infty} \sum_{i=0}^{p^n-1}g(i)\mu_{P, -1}(i+p^n \Bbb Z_p)  \notag\\
&=&  \lim_{n \rightarrow \infty} \sum_{i=0}^{p^n-1}g(i)(-1)^i i^k \\
&=& \int_{\Bbb Z_p} g(x) x^k d \mu_{-1} (x). \notag
\end{eqnarray}
Therefore we have the following theorem.

\begin{theorem} Let $P(x) \in \mathbb C_p [[x]]$ be an arbitrary polynomial.
Then $\mu_{P, -1}$ is a strongly fermionic $p$-adic invariant
measure and for all $x \in \Bbb Z_p$,
$$f_{\mu_{P, -1}}=(-1)^x P(x).$$
Furthermore, for any $g \in UD(\Bbb Z_p, \Bbb C_p)$,
$$\int_{\Bbb Z_p} g(x) d \mu_{P, -1} (x)=\int_{\Bbb Z_p} g(x) P(x) d \mu_{-1} (x),$$
where the second integral is fermionic $p$-adic invariant integral
on $\Bbb Z_p$.
\end{theorem}

Let $f(x)=\sum_{n=0}^{\infty} a_n \binom{x}{n}$ be the Mahler
expansion of uniformly differential function of $f$, where
$\binom{x}{n}$ is the binomial coefficient. Then we note that
$\lim_{n \rightarrow \infty} n |a_n | =0.$ Let $f_m
(x)=\sum_{i=0}^{m} a_i \binom{x}{i} \in \Bbb C_p [[x]].$ Then
\begin{eqnarray}||f- f_m ||_\infty \le \underset{n \ge m}{\sup} \,n|a_n|.\end{eqnarray}
Writing $f=f_m+ f-f_m$, we easily see that
\begin{eqnarray*}
& & |\mu_{f, -1}(a+p^n \Bbb Z_p )-\mu_{f, -1}(a+p^{n+1} \Bbb Z_p )| \\
& &\le \max \{ |\mu_{f_{m, -1}}(a+p^n \Bbb Z_p)-\mu_{f_{m,
-1}}(a+p^{n+1} \Bbb Z_p)|, \\ & & \qquad \quad \, \, |\mu_{f-f_{m,
-1}}(a+p^n \Bbb Z_p)-\mu_{f-f_{m, -1}}(a+p^{n+1} \Bbb Z_p)| \}.
\end{eqnarray*}

By Theorem 1, we note that for some constant $C_1$,
\begin{eqnarray} | \mu_{f-f_m, -1} (a+p^n \Bbb Z_p)| \le || f-f_m ||_\infty \le C_1 p^{-n}.\end{eqnarray}
Also, for $m \gg 0$, it follows that $|| f||_\infty = ||f_m
||_\infty $, and so
\begin{eqnarray}|\mu_{f_{m, -1}}(a+p^n \Bbb Z_p) - \mu_{f_{m, -1}}(a+p^{n+1} \Bbb Z_p)| \le C_2 p^{-n},\end{eqnarray}
where $C_2$ is some positive constant.

From (14), we note that
\begin{eqnarray*}
& & |f(a)-\mu_{f, -1}(a+p^{n} \Bbb Z_p )| \\
& &\le \max \{ |f(a)-f_m (a)|, \, |f_m (a)-\mu_{f_m, -1} (a+p^{n}
\Bbb Z_p)|, \, |\mu_{f-f_m, -1} (a+p^{n} \Bbb Z_p) | \} \\
& &\le \max \{ |f(a)-f_m (a)|, \, |f_m (a)-\mu_{f_m, -1} (a+p^{n}
\Bbb Z_p)|, \,  ||f-f_m ||_\infty \}.
\end{eqnarray*}
If we fix $\epsilon >0$, and fix $m$ such that $||f-f_m ||_\infty
\le \epsilon$, then for $n  \gg 0$, we have
$$|f(a)-\mu_{f, -1}(a+p^n \Bbb Z_p)| \le \epsilon.$$
Hence we see that
\begin{eqnarray}f_{\mu_{f,-1}}(a) = \lim_{n \rightarrow \infty} \mu_{f, -1} (a+p^n \Bbb Z_p)=(-1)^a f(a). \end{eqnarray}
Let $m$ be the sufficiently large number such that $||f-f_m
||_\infty \le p^{-n} $. Then we have
\begin{eqnarray*}
\mu_{f, -1} (a+p^n \Bbb Z_p)&=&\mu_{f_m , -1} (a+p^n \Bbb
Z_p)+\mu_{f-f_m , -1} (a+p^n \Bbb Z_p)\\
& \equiv & \mu_{f_m , -1} (a+p^n \Bbb Z_p)\equiv (-1)^a f(a)\,E_0 \,
(\text{mod} \, p^n).
\end{eqnarray*}

Let $g \in UD(\Bbb Z_p, \Bbb C_p)$. Then
$$\int_{\Bbb Z_p} g(x) d \mu_{f, -1} (x)=\int_{\Bbb Z_p} f(x) g(x) d \mu_{-1} (x) .$$

Let $f$ be the function from $UD(\Bbb Z_p, \Bbb C_p)$ to Lip$(\Bbb
Z_p, \Bbb C_p)$. We know that $\mu_{-1}$ is a strongly $p$-adic
invariant measure on $\Bbb Z_p$, and
\begin{eqnarray}|f_{\mu_{-1}}(a)-\mu_{-1} (a+p^n \Bbb Z_p ) | \le C_4 p^{-n} \end{eqnarray}
for any positive integer $n$ and some constant $C_4$.

If $\mu_{1, -1}$ is the associated strongly fermionic $p$-adic
invariant measures on $\Bbb Z_p$, then we have
\begin{eqnarray}|\mu_{1, -1} (a+p^n \Bbb Z_p)-f_{\mu_{- 1}}(a)| \le C_5 p^{-n} \end{eqnarray}
for $n\gg0$. Thus, we note that for sufficiently large $n$,
\begin{eqnarray}
& & |\mu_{-1} (a+p^n \Bbb Z_p)-\mu_{1, - 1}(a+p^n \Bbb Z_p)| \notag \\
& & \quad \le |\mu_{-1} (a+p^n \Bbb
Z_p)-f_{\mu_{-1}}(a)|+|f_{\mu_{-1}}(a)- \mu_{1, - 1}(a+p^n \Bbb
Z_p)| \le K, \end{eqnarray} where $K$ is a constant. Therefore
$\mu_{-1} - \mu_{1, -1}$ is a measure on $\Bbb Z_p$. Hence we obtain
the following theorem.

\begin{theorem} Let $\mu_{-1}$ be a strongly fermionic $p$-adic invariant
on $\Bbb Z_p$, and assume that the fermionic Radon-Nikodym
derivative $f_{\mu_{-1}}$ on $\Bbb Z_p$ is uniformly differentiable
function. Suppose that $\mu_{1, -1}$ is the strongly fermionic
$p$-adic invariant measure associated to $f_{\mu_{-1}}$. Then there
exists a measure $\mu_{2, -1}$ on $\Bbb Z_p$ such that
$$\mu_{-1} = \mu_{1,-1}+\mu_{2,-1}.$$
\end{theorem}

\begin{remark} Theorem 2 seems to be the $p$-adic analogue of
Lebesque decomposition with respect to strongly fermionic $p$-adic
invariant measure on $\Bbb Z_p$, ( see [9, 10]).
\end{remark}

\end{document}